\documentclass[12pt]{amsart}
\usepackage{amsmath,amssymb,amscd,verbatim}
\usepackage{amsfonts}

\newcommand{\mc}{\mathcal}
\newcommand{\sub}{\subseteq}
\newcommand{\ol}{\overline}

\newcommand{\sm}{\setminus}
\DeclareMathOperator{\rad}{rad}

\newtheorem{theorem}{Theorem}[section]
\newtheorem{lemma}[theorem]{Lemma}
\newtheorem{proposition}[theorem]{Proposition}

\theoremstyle{definition}
\newtheorem{definition}[theorem]{Definition}
\newtheorem{remark}[theorem]{Remark}
\newtheorem{example}[theorem]{Example}

\newcommand\mpb{\medskip}
\newcommand\bpb{\bigskip}

\newcommand\tpp{\vskip5pt}

\DeclareMathOperator{\Spec} {\mbox{\rm Spec}}
\DeclareMathOperator{\Max} {\mbox{\rm Max}}

   \newcommand{\balf}
   {\renewcommand{\theenumi}{(\alph{enumi})}
   \renewcommand{\labelenumi}{\theenumi}
                        \begin{enumerate}}
  \newcommand{\ealf}   {\end{enumerate}
                        \renewcommand{\theenumi}{\arabic{enumi}}
                        \renewcommand{\labelenumi}{\theenumi.}}
  \newcommand{\bara}   {\renewcommand{\theenumi}{(\arabic{enumi})}
                        \renewcommand{\labelenumi}{\theenumi}
                        \begin{enumerate} }
  \newcommand{\eara}   {\end{enumerate}
                        \renewcommand{\theenumi}{\arabic{enumi}}
                        \renewcommand{\labelenumi}{\theenumi.}}

   \newcommand{\brom} {\renewcommand{\theenumi}{
(\roman{enumi})} \renewcommand{\labelenumi}{\theenumi}
                        \begin{enumerate} }

  \newcommand{\erom}   {\end{enumerate}
                        \renewcommand{\theenumi}{\arabic{enumi}}
                        \renewcommand{\labelenumi}{\theenumi.}}

\begin{document}
\title{Factoring Ideals in Pr\"ufer Domains}

\author{Marco Fontana}
\address{Dipartimento di Matematica\\ Universit\`{a} degli Studi Roma Tre\\
    Largo San L. Murialdo, 1\\
    00146 Roma, Italy}
\email{fontana@mat.uniroma3.it}
\author{Evan Houston}
\address[Evan Houston and Thomas Lucas]{Department of Mathematics\\
    University of North Carolina at Charlotte\\
    Charlotte, NC 28223 U.S.A.}
\email[Evan Houston]{eghousto@email.uncc.edu}
\author{Thomas Lucas}
\email[Thomas Lucas]{tglucas@email.uncc.edu}
\thanks{The first author was supported by MIUR, under Grant
PRIN 2005-015278, and the second author was supported by a
visiting grant from GNSAGA of INdAM (Istituto Nazionale di Alta
Matematica).}

\date{}

\begin{abstract} We show that in certain Pr\"ufer domains, each
nonzero ideal $I$ can be factored as $I=I^v \Pi$, where $I^v$ is
the divisorial closure of $I$ and $\Pi$ is a product of maximal
ideals.  This is always possible when the Pr\"ufer domain is
$h$-local, and in this case such factorizations have certain
uniqueness properties.  This leads to new characterizations of the
$h$-local property in Pr\"ufer domains. We also explore
consequences of these factorizations and give illustrative
examples.
\end{abstract}

\maketitle

\baselineskip16pt

\mpb


Let $R$ be a Pr\"ufer domain.  Recall that $R$ has finite
character if each nonzero element of $R$ is contained in only
finitely many maximal ideals of $R$ and that $R$ is $h$-local if
it has finite character and each nonzero prime ideal of $R$ is
contained in a unique maximal ideal of $R$.  It follows from
\cite[Theorem 4.12]{bhlp} that if $R$ is $h$-local, then each
nonzero ideal $I$ of $R$ factors as $I=I^v\Pi$, where $I^v$
denotes the divisorial closure of $I$ and $\Pi$ is a product of
maximal ideals.  Part of the first section of this work may be
viewed as an elaboration of this result.  We observe that, for a
nonzero ideal $I$ of an $h$-local Pr\"ufer domain, we have
$I=I^vM_1 \cdots M_n$, where the $M_i$ are precisely the
nondivisorial maximal ideals $M$ of $R$ which contain $I$ and for
which $IR_M$ remains nondivisorial in $R_M$ (and where we take the
empty product of maximal ideals to be $R$ itself); moreover, this
factorization is unique in the sense that no $M_i$ can be deleted.
On the other hand, we show that in certain almost Dedekind
domains, one can have a weaker factorization property: each
nonzero ideal $I$ factors as $I=I^v\Pi$, where $\Pi$ is a product
of (not necessarily distinct) maximal ideals. We show
(Proposition~\ref{p:nonmaxdiv}) that in a Pr\"ufer domain with
this weak factorization property each nonmaximal prime ideal is
divisorial, each branched nonmaximal prime ideal is the radical of
a finitely generated ideal, and each branched idempotent maximal
ideal is sharp. (Relevant definitions are reviewed in the sequel.)
If, in addition to possessing the weak factorization property, the
Pr\"ufer domain $R$ has finite character, then $R$ is $h$-local
(Theorem~\ref{t:converse}). Moreover, a Pr\"ufer domain is
$h$-local if and only if it has the strong factorization property
(Theorem \ref{t:bigthm}). Another interesting property of
$h$-local Pr\"ufer domains is that a nonzero ideal of such a
domain is divisorial if and only if it is locally divisorial (at
maximal ideals). In fact, we show in Theorem~\ref{t:bigthm} that a
Pr\"ufer domain with this property is $h$-local.

In Section 2, we explore in $h$-local Pr\"ufer domains how a given
factorization of an ideal $I$ affects that of $\rad I$ and
$II^{-1}$ and how factorizations of ideals $I$ and $J$ affect
those of such related ideals as $IJ$, $I \cap J$, and $I+J$.

Section 3 is devoted to examples.  As has already been mentioned,
it is possible for an almost Dedekind domain to possess the weak
factorization property; in Example~\ref{e:infnondiv} we show that
this can happen even in an almost Dedekind domain with infinitely
many nondivisorial maximal ideals. While in a Pr\"ufer domain with
the strong factorization property, the sum of two divisorial
ideals must be again divisorial, we show in
Example~\ref{e:sumnotdiv} that an almost Dedekind domain may have
the weak factorization property and still possess divisorial
ideals $I$ and $J$ with $I+J$ not divisorial. We also give an
example (Example~\ref{e:nofac}) of a one-dimensional Bezout domain
$R$ which does not have the weak factorization property, and we
observe that in this example, there is a divisorial ideal $J$ and
a maximal ideal $M$ with $JR_M$ not divisorial. \bigskip

\emph{Acknowledgement.} The authors would like to thank Bruce
Olberding and the referee, whose many helpful comments greatly
improved this paper.


\section{The strong and weak factorization properties}

We begin by recalling some facts which we shall use frequently and
without further reference.  Let $V$ be a valuation domain with
maximal ideal $M$.  If $M$ is divisorial, then $M$ is principal
and every nonzero ideal of $V$ is divisorial by \cite[Lemma
5.2]{h}. On the other hand, if $M$ is not divisorial, then by
\cite[Lemma 4.2]{bg} a nonzero ideal $I$ of $V$ is nondivisorial
if and only if $I=xM$ for some element $x \in V$.

\begin{theorem} \label{t:hlocfac}
Let $R$ be an $h$-local Pr\"ufer domain. Then
\begin{enumerate}
\item[(1)] For each divisorial ideal $I$ of $R$, if $M\supseteq I$ with $M$ a
non-divisorial maximal ideal of $R$, then $IR_M$ is divisorial in
$R_M$, and $IR_M$ is properly contained in $MR_M$.
\item[(2)] For each nonzero nondivisorial ideal $I$ of $R$, $I$
factors as a product $BM_1M_2\cdots M_n$ where  $B$ is a
divisorial ideal and the $M_i$ are  distinct non-divisorial
maximal ideals of $R$ that contain $I$ for which   $IR_{M_i}$ is
not a divisorial ideal of $R_{M_i}$. Moreover, this factorization
is unique in the sense that $B = I^v$ and the $M_i$ include all
maximal ideals that contain $I$ where $IR_{M_i}$ is not
divisorial.
\end{enumerate}
\end{theorem}

\begin{proof} Let $A$ be a nonzero ideal of $R$. Since $R$ is $h$-local,
$(AR_M)^{-1}=A^{-1}R_M=(A^vR_M)^{-1}$ for each maximal ideal $M$
(\cite[Lemma 2.3]{bs} or \cite[Theorem 3.10]{o}). Moreover,
$A^vR_M=(AR_M)^v$. In particular, if $M$ is nondivisorial, then
$(MR_M)^v=M^vR_M=R_M$, so that $MR_M$ is not divisorial, while if
$I$ is divisorial, then $IR_M$ is also divisorial.  This proves
(1).

If $AR_M$ is not divisorial, then it must be of the form $xMR_M$
for some $x\in R$. In this case, we have $A^vR_M=(AR_M)^v=xR_M$
and $AR_M=A^vMR_M$.

Now let $I$ be a nonzero nondivisorial ideal of $R$. Let $M_1$,
$M_2$, \dots , $M_n$ be the nondivisorial maximal ideals that
contain $I$ where $IR_{M_i}$ is not divisorial. (It will follow
from the rest of the proof that $n >0$, but for the moment we take
the empty product to be $R$.) Consider the ideal
$J=I^vM_1M_2\cdots M_n$. For each $M_i$, it is clear that
$JR_{M_i}=IR_{M_i}$ from the argument above. Let $M$ be any other
maximal ideal. If $M$ does not contain $I$, then $JR_M=R_M=IR_M$.
On the other hand if $M$ contains $I$, we must have that
$(IR_M)^v=IR_M$. As $I^vR_M=(IR_M)^v$, we obtain
$IR_M=I^vR_M=JR_M$. Hence $I=J$.

Now suppose $I=BN_1N_2\cdots N_m$ with $B$ divisorial and the
$N_i$ distinct members of $\{M_1,M_2,\dots,M_n\}$. Since for each
$i$, $BR_{M_i}$ is divisorial (but perhaps trivial) and $IR_{M_i}$
is not, checking locally at $M_i$ shows that some $N_j$ must equal
$M_i$. Hence $m=n$ and each $M_i$ is needed in the factorization.
Rewriting, we have $I=BM_1M_2\cdots M_n$. Thus, since the $M_i$
are nondivisorial (and since for a maximal ideal $M$, we have $M$
nondivisorial if and only if $M^{-1}=M^v=R$), $I^v=(BM_1 \cdots
M_n)^v=B^v=B$.
\end{proof}

\mpb

\begin{definition} \label{d:sfp} A Pr\"ufer domain $R$ has the
\emph{strong factorization property} if for each nonzero ideal $I$
of $R$, we have (1) $I=I^vM_1 \cdots M_n$ where $M_1, \ldots, M_n$
are precisely the nondivisorial maximal ideals of $R$ which
contain $I$ for which $IR_M$ is nondivisorial and (2) this
factorization is unique in the sense that no $M_i$ can be omitted.
\end{definition}

\begin{remark} \label{r:emptyprod} In Definition~\ref{d:sfp}, we
take the empty product of maximal ideals to be $R$; thus, if
$IR_M$ is divisorial for each maximal ideal $M$, then $I=I^v$
(that is, $I$ is divisorial).
\end{remark}

Thus, according to Theorem~\ref{t:hlocfac}, $h$-local Pr\"ufer
domains possess the strong factorization property.  In
Theorem~\ref{t:bigthm} below, we show that the converse holds.

\begin{remark} \label{r:sfp} Let $I$ be a nonzero ideal of the Pr\"ufer domain
$R$, denote by $\Max(R,I)$ the set of maximal ideals of $R$
containing $I$, and set
           $$
           \begin{array}{rl}
               \mathfrak{M}(I) := &\hskip -5pt  \{M  \in \Max(R, I)  \mid M^v =
               R\,,\; IR_{M}\neq (IR_{M})^v\}\\
               \mathfrak{M'}(I) := &\hskip -5pt
                \{M'  \in \Max(R, I)  \mid {M'}^v =
               R\,,\; IR_{M'} =(IR_{M'})^v\}\\
                 \mathfrak{N}(I) := &\hskip -5pt
                \{N \in \Max(R, I)  \mid N = {N}^v  \}\,.
    \end{array}
             $$

             Then Definition~\ref{d:sfp} requires that $\mathfrak{M}(I)$ be
             finite (possibly empty),
             that $I = I^v \prod_{M \in \mathfrak{M}(I)} M$, and
             that this factorization be irredundant.  We say
             nothing about the possible finiteness of
             $\mathfrak{M'}(I)$ or $\mathfrak{N}(I)$.  It is also
             possible that $I$ could have a different
             factorization involving some of the maximal ideals
             in $\mathfrak{M'}(I) \cup \mathfrak{N}(I)$.
             For example, let $(V,M)$ be a valuation domain
             containing a non-principal divisorial ideal $I$.
             Then $\mathfrak{M}(I)$ is empty, and the
             factorization of $I$ is just $I=I^v$.  However,
             since $I$ is not principal, we also have $I=IM
             (=I^vM)$.  (The fact that $I$ not principal implies that
             $I=IM$ is probably well known, but here is a proof: Begin
             with an element $x \in I$.  Since $I$ is not principal, we
             may then choose $y \in I \sm Vx$ so that $x/y \in M$ and
             $x = y(x/y) \in IM$.) By constructing $V$ appropriately, we may
             have $M$ divisorial or not, that is, $\mathfrak{N}(I)= \{M\}$
             or  $\mathfrak{M'}(I) = \{M\}$.
\end{remark}

\begin{remark} \label{r:sfp2} Using the notation in
Remark~\ref{r:sfp} and following the proof of \cite[Theorem
4.12]{bhlp}, we have for any nonzero ideal $I$ in an $h$-local
Pr\"ufer domain a decomposition of $I^v$ as follows.  Set
$I'=\bigcap_{M' \in \mathfrak{M'}(I)} (IR_{M'} \cap R)$ and $J_I=
\bigcap_{N \in \mathfrak{N}(I)} (IR_N \cap R)$.  For each $M \in
\mathfrak{M}(I)$, there is an invertible ideal $L$ of $R$ with
$IR_M \cap R=LM$; set $L_I$ equal to the product of these $L$'s.
Then $I^v=L_II'J_I$. \end{remark}

We now introduce our second factorization property.

\begin{definition} \label{d:wfp} A Pr\"ufer domain $R$ has the
\emph{weak factorization property} if each nonzero ideal $I$ can
be written as $I=I^v \Pi$, where $\Pi$ is a (finite) product of
(not necessarily distinct) maximal ideals (and where, again, the
empty product of maximal ideals is taken to be $R$).
\end{definition}

Before stating our next few results, we need some terminology.
Recall that a domain $R$ satisfies the \emph{trace property} if,
for each nonzero ideal $I$ of $R$, we have that $II^{-1}$ is equal
either to $R$ or to a prime ideal of $R$.  The domain $R$
satisfies the \emph{radical trace property} if each nonzero ideal
$I$ of $R$ satisfies $II^{-1}=R$ or $II^{-1}=\rad(II^{-1})$.
Finally, $R$ satisfies the \emph{weak trace property for primary
ideals} if, for each nonzero, nonmaximal prime ideal $P$ and each
$P$-primary ideal $Q$, we have $QQ^{-1}=P$.  For information about
the trace and radical trace properties, the reader is referred to
\cite{fhp} and \cite{lu}.  Now recall from \cite{g} that a domain
$R$ is said to be a \emph{$\#$-domain} if $\bigcap_{M \in \mathcal
M} R_M \ne \bigcap_{N \in \mathcal N} R_N$ for each pair of
distinct nonempty subsets $\mathcal M$ and $\mathcal N$ of the set
of maximal ideals of $R$, equivalently, if for each maximal ideal
$M$ of $R$, $R_M$ does not contain $\bigcap R_N$, where the
intersection is taken over those maximal ideal $N$ with $N \ne M$.
This was extended to focus on a single maximal ideal in \cite{ll}:
a maximal ideal is \emph{sharp} if $R_M$ does not contain
$\bigcap_{N \ne M} R_N$. By \cite[Corollary 2]{gh} a maximal ideal
$M$ of a Pr\"ufer domain $R$ is sharp if and only if there is a
finitely generated ideal of $R$ which is contained in $M$ and no
other maximal ideal of $R$. Finally, a domain $R$ is a
\emph{$\#\#$-domain} if each overring of $R$ is a $\#$-domain (see
\cite{gh}).

\begin{proposition} \label{p:nonmaxdiv} Let $R$ be a Pr\"ufer domain with the
weak factorization property.  Then
\begin{enumerate}
\item[(1)] each ideal which is primary to a nonmaximal ideal of $R$ is
divisorial (in particular, each nonmaximal prime is divisorial),
\item[(2)] if $M$ is an idempotent maximal ideal of $R$ and $I$ is a
nondivisorial $M$-primary ideal, then $I=I^vM$,
\item[(3)] each branched maximal idempotent ideal of $R$ is sharp,
\item[(4)] $R$ has the weak trace property for primary ideals, and
\item[(5)] each branched nonmaximal prime ideal of $R$ is the
radical of a finitely generated ideal.
\end{enumerate}
\end{proposition}

\begin{proof} (1)  Let $Q$ be a $P$-primary ideal of $R$ with $P$
nonmaximal. Write $Q=Q^v \Pi$, where $\Pi$ is a product of maximal
ideals.  Then $\Pi \nsubseteq P$, whence $Q^v \sub Q$, and so $Q$
is divisorial.

(2) This is clear.

(3) Let $M$ be a branched idempotent maximal ideal of $R$. Since
$M$ is branched, there is an $M$-primary ideal $I$ with $I \ne M$.
If $I^v \nsubseteq M$, then $I^v =R$, and $I=I^vM$ by (2).  But
this yields $I=M$, a contradiction.  Hence $I^v \sub M$, and $M$
is sharp by \cite[Proposition 2.2]{o}.

(4)  Let $Q$ be a proper $P$-primary ideal with $P$ not maximal.
Then $Q$ is divisorial by (1).  We shall show that $QQ^{-1}=P$. By
\cite[Corollary 3.1.8 and Theorem 3.1.2]{fhp} $P^{-1} = \bigcap
R_M \cap R_P$, where the intersection $\bigcap R_M$ is taken over
those maximal ideals which do not contain $P$. For $x \in Q^{-1}$,
we have $(R:_R x) \nsubseteq M$, since $Q \nsubseteq M$; thus $x
\in R_M$.  Hence $Q^{-1} \sub \bigcap R_M$.  The same argument
shows that $Q^{-1} \sub \Omega(P) := \bigcap R_N$, where $N$
ranges over the prime ideals of $R$ which do not contain $P$.

For $y \in P^{-1}$, we have $y \in R_P$, whence $ay \in R$ for
some $a \notin P$.  Then $ayQ \sub Q$ yields $yQ \sub Q$ (since it
is clear that $yQ \sub R$). Thus $P^{-1}Q \sub Q$. Therefore,
$(QQ^{-1})^{-1}=(QQ^{-1}:QQ^{-1}) \supseteq P^{-1}$, and we have
$QQ^{-1} \sub P^v=P$ by (1).  We also have that $P^{-1} \sub
(QQ^{-1})^{-1} \sub Q^{-1} \sub \Omega(P)$ with $Q^{-1}$ a ring.
Since $R$ is a Pr\"ufer domain, this yields $(QQ^{-1})^{-1} =
P^{-1}$ \cite[Theorem 3.3.7]{fhp}, whence $(QQ^{-1})^v = P$ (again
by (1)). If $QQ^{-1}$ is not divisorial, then $QQ^{-1}=(QQ^{-1})^v
\Pi = P \Pi$, for some product $\Pi$ of maximal ideals each of
which necessarily contains $P$ (since each contains $Q$).  A
routine local check then shows that $P \Pi = P$, so that $QQ^{-1}
= P$, as desired.

(5) This follows from (1) and \cite[Proposition 2.9]{fhl}.
\end{proof}

Next, we give some consequences of the strong factorization
property.

\begin{theorem} \label{t:sfp2} Let $R$ be a Pr\"ufer domain
with the strong factorization property. Then
\begin{enumerate}
\item If $I$ is a nonzero ideal of $R$, then $I$ is divisorial
if and only if $IR_M$ is divisorial for each maximal ideal $M$ of
$R$.
\item If $M$ is a maximal ideal of $R$ and $A$ is a
divisorial ideal of $R_M$, then $A \cap R$ is divisorial in $R$.
\item For each maximal ideal $M$, if $M$ is not divisorial, then $MR_M$
is not divisorial. Thus the non-divisorial maximal ideals are those that are
idempotent.
\item For each nonzero ideal $I$ of $R$ and each maximal ideal $M$
of $R$, we have $(IR_M)^v=I^vR_M$.
\item If $I$ is an ideal contained in no nondivisorial maximal ideals, then $I$ is divisorial.
\item For each fractional ideal $F$, $F=F^vM_1M_2\cdots M_n$ where the
$M_i$ are the maximal ideals that contain some particular fixed nonzero
principal multiple $xF$ of $F$ with $xFR_{M_i}$ not divisorial. Moreover,
the factorization is unique.
\item   If $R$ has finite character, and $I$ is an ideal for which $IR_M$ is divisorial
only in the trivial case $IR_M=R_M$, then $I^v$ is invertible.

\end{enumerate}\end{theorem}

\begin{proof}

(1)  Let $I$ be a nonzero ideal of $R$, and let $M$ be a maximal
ideal.  Suppose that $I$ is divisorial.  If $M$ is nondivisorial,
then $IR_M$ is divisorial by Definition~\ref{d:sfp}.  If $M$ is
divisorial, then it is invertible; hence $MR_M$ is principal, and
every ideal of $R_M$ is divisorial.  For the converse, see
Remark~\ref{r:emptyprod}.

(2) Let $M$ be maximal, and let $A$ be a divisorial ideal of
$R_M$. Set $I=A \cap R$, and write $I=I^vM_1 \cdots M_n$ as in
Definition~\ref{d:sfp}.  Since $IR_M=A$ is divisorial, $M \notin
\{M_1, \cdots M_n\}$.  We then have $IR_M = I^vM_1 \cdots M_n R_M
= I^vR_M$.  Hence $I^v \sub I^vR_M \cap R = IR_M \cap R = I$, as
desired.

(3)  From (1) if $M$ is a nondivisorial maximal ideal, then $MR_M$
is also nondivisorial and hence idempotent.  Since idempotence is
a local property, $M$ is itself idempotent.

(4) Let the factorization of $I$ be $I=I^vM_1 \cdots M_n$, and let
$M$ be a maximal ideal of $R$.  If $M \notin \{M_i\}$, then
$(IR_M)^v=(I^vM_1 \cdots M_nR_M)^v=(I^vR_M)^v=I^vR_M$, with the
last equality following from (1). If $M=M_i$ for some $i$, then
$(IR_M)^v=(I^vM_1 \cdots M_nR_M)^v=(I^vMR_M)^v=(I^vR_M)^v=I^vR_M$,
with the penultimate equality following from (2) and the last
equality following from (1).

(5)  This is immediate from the definition.

(6)  Let $F$ be a fractional ideal and let $x\in R \sm \{0\}$ be
such that $xF\subseteq R$. Then we can factor $xF$ uniquely as
$(xF)^vM_1M_2\cdots M_n$ where the $M_i$ are the nondivisorial
maximal ideals that contain $xF$ where $xFR_{M_i}$ is not
divisorial. Of course, $(xF)^v=xF^v$, so we can cancel the $x$ to
obtain $F=F^vM_1M_2\cdots M_n$. For any nonzero $y\in (R:F)$, we
obtain a (possibly different) factorization $F=F^vN_1N_2\cdots
N_k$ where the $N_j$ are such that $yFR_{N_j}$ is not divisorial.
If these two factorizations were actually different, we would have
two distinct factorizations of $xyF$, one as $xyF^vM_1M_2\cdots
M_n$ and the other as $xyF^vN_1N_2\cdots N_k$. Thus we must have a
unique factorization for $F$.

(7)  Let $I$ be as indicated.  Then $I=I^vM_1 \cdots M_n$, where
the $M_i$ are precisely the maximal ideals which contain $I$.  For
each $i$, $IR_M$ not divisorial yields an element $x_i \in I^v$
with $IR_{M_i}=x_iM_iR_{M_i}$ and $I^vR_{M_i}=x_iR_{M_i}$.  Let
$A=(x_1,x_2,\dots,x_n)$. At most  finitely many maximal ideals
contain $A$,  say $N_1$, $N_2$, \dots, $N_k$. For those $N_j$ that
are not among the $M_i$s, we may choose an element $y_j\in I \sm
N_j$. Let $B$ be the ideal generated by $A$ and the $y_j$.
Obviously, $B\subseteq I^v$.  Now consider the ideal
$J=BM_1M_2\cdots M_n$ and let $M$ be a maximal ideal of $R$. If
$M=M_i$ for some $i$, then $JR_{M_i}=BM_iR_{M_i}$. Thus  $JR_{M_i}
= I^vM_iR_{M_i} = IR_{M_i}$ since $x_iR_{M_i}=I^vR_{M_i}$ and
$x_i\in B\subseteq I^v$. If $M$ is a maximal ideal not among the
$M_i$, then $B \nsubseteq M$, and we have $JR_M = BR_M = R_M =
IR_M$ since no other maximal ideals contain $B$. Hence $J=I$. As
$B$ is divisorial and factorizations are unique, we must have
$B=I^v$. Therefore, $I^v$ is invertible.
\end{proof}

\mpb

We observe that, in view of Theorem~\ref{t:bigthm} below, part (5)
of Proposition~\ref{t:sfp2} is \cite[Proposition 6.5 (a)]{hhp} and
part (7) may be viewed as a generalization of \cite[Proposition
6.5 (b)]{hhp}. \medskip

We need a couple of general results before proving that statement
(1) in Theorem~\ref{t:sfp2} is equivalent to the $h$-local
property.  Our next lemma provides a way to prove statement (2) of
Theorem~\ref{t:sfp2} using only the assumption that each locally
divisorial of the Pr\"ufer domain $R$ is divisorial.

\begin{lemma} \label{l:mnlocdiv} Let  $I$ be a nonzero ideal of
a Pr\"ufer domain $R$ and let $M$ a maximal ideal that contains
$I$.  For $J=IR_M\bigcap R$, $JR_N$ is a divisorial ideal of $R_N$
for each maximal ideal $N\ne M$.
\end{lemma}
\begin{proof} Let $N$ be a maximal ideal of $R$ with $N \ne M$. Then
$JR_N = (IR_M\bigcap R)R_N = IR_MR_N\bigcap R_N = IR_P\bigcap R_N$
where $P$ is the largest prime contained in $M \cap N$.  If $JR_N$
is not divisorial, then $JR_N=xNR_N$ for some $x \in R$.  This
yields $JR_P=xR_P$, and we then have $x \in JR_P \bigcap R_N =
IR_P \bigcap R_N = JR_N = xNR_N$, a contradiction.  Hence $JR_N$
is divisorial.
\end{proof}

\begin{theorem} \label{t:fgcondiv}  Let   $R$  be a Pr\"ufer
domain and let $P$ be a nonzero nonmaximal prime  that is the
radical of a finitely generated ideal. If $I$ is a  finitely
generated ideal whose radical is $P$ and  $M$ is a maximal ideal
that contains $P$, then the ideal  $J=IR_M\bigcap R$  is
divisorial if and only if $M$ is the only maximal ideal that
contains $P$. \end{theorem}

\begin{proof} Let $J=IR_M\bigcap R$ where $M$ is a maximal ideal
that contains $P$. It is clear that if $M$ is the only maximal
ideal that contains $P$, then $J^v=J=I$.

For the remainder of the proof, we assume that $M$ is not the only
maximal ideal that contains $P$. Denote by $P'$ the largest prime
ideal contained in all the maximal ideals which contain $I$. Then
$P'$ is properly contained in $M$. We shall show that
$J^{-1}=P'I^{-1}$.

We check the inclusion $P'I^{-1} \sub J^{-1}$ locally.  At $M$ we
have $I^{-1}P'JR_M = I^{-1}P'IR_M \sub R_M$. For $N \in \Max(R,I)
\setminus \{M\}$, we have $I^{-1}P'JR_N = I^{-1}P'(IR_{P'} \cap
R_N) \sub I^{-1}P'IR_{P'} = I^{-1}IP'R_N \sub R_N$.  Finally, for
$L \notin \Max(R,I)$, we have $I^{-1}P'JR_L = I^{-1}R_L
=(IR_L)^{-1} = R_L$.  Thus $P'I^{-1} \sub J^{-1}$.

For the reverse inclusion, let $t \in J^{-1}$.  Choose any $N \in
\Max(R,I) \setminus \{M\}$, and then choose $a \in NR_N \setminus
P'R_N$.  Then $a^{-1}I \sub IR_{P'} \cap R_N=JR_N$. Hence
$ta^{-1}I \sub tJR_N \sub R_N$, yielding $tI \sub aR_N$.  It
follows that $tI \sub P'R_N \cap R = P'$.  Thus $J^{-1}I \sub P'$,
and we have $J^{-1} \sub I^{-1}P'$, as desired.

Finally, we show that $J$ is not divisorial.  Suppose, on the
contrary, that $J = J^v = IP'^{-1}$.  Then $I^{-1}J = P'^{-1}$.
Now choose $m \in M \setminus P'$, and then choose $u \in
(I,m)^{-1} \setminus R_M$.  Then $(R:_R u) \nsubseteq P'$ and
$(R:_R u) \nsubseteq L$ for each maximal ideal $L$ with $L \notin
\Max(R,I)$.  It follows that $u \in R_{P'} \cap \left(\bigcap
\{R_L \mid L \notin \Max(R,P')\}\right) = P'^{-1}$ \cite[Theorem
3.1.2 and Corollary 3.1.8]{fhp}.  Hence $u \in P'^{-1}R_M =
I^{-1}JR_M=R_M$, a contradiction.
\end{proof}

\begin{lemma} \label{l:locdivrtp} Let $R$ be a Pr\"ufer domain.
If $R$ has the property that an ideal $I$ of $R$ is divisorial
whenever $IR_M$ is divisorial for each maximal ideal $M$, then $R$
has the radical trace property.
\end{lemma}
\begin{proof} Assume that $R$ has the property that each locally
divisorial ideal is divisorial. By \cite[Theorem 23]{lu}, to show
that $R$ has the radical trace property, it suffices to show if
$Q$ is a $P$-primary ideal such that $Q^{-1}$ is a ring, then
$Q=P$. To this end, let $Q$ be a proper $P$-primary ideal. Since
$R$ is integrally closed, $Q^{-1}$ is a ring if and only if
$Q^{-1}=P^{-1}=(P:P)$ \cite[Proposition 3.1.16]{fhp}.

If $P$ is not maximal, then $QR_M$ is divisorial for each maximal
ideal $M$ (see the argument that $JR_N$ is divisorial in
Lemma~\ref{l:mnlocdiv} above). Hence $Q$ is divisorial and we have
$P^{-1}\subsetneq Q^{-1}$. Thus $Q^{-1}$ is not a ring.

If $P$ is maximal and $Q$ is divisorial, then we again have that
$Q^{-1}$ is not a ring. The only other case is when $QR_P=xPR_P$
with $P$ idempotent and $x$ some nonzero element of $P$. Then
$Q'=xR_P\bigcap R$ is a proper $P$-primary ideal which is
divisorial since it is divisorial in each $R_N$. Hence we have
$P^{-1}\subsetneq Q'^{-1}\subseteq Q^{-1}$, and again $Q^{-1}$ is
not a ring. \end{proof}

\begin{theorem} \label{t:bigthm} The following are equivalent for
a Pr\"ufer domain $R$.
\begin{enumerate}
\item $R$ is $h$-local.
\item $R$ has the strong factorization property.
\item For each nonzero ideal $I$ of $R$, $I$ is divisorial if and
only if $IR_M$ is divisorial in $R_M$ for each maximal ideal $M$
of $R$.
\item For each nonzero ideal $I$ of $R$, if $IR_M$ is divisorial
for each maximal ideal $M$, then $I$ is divisorial.
\end{enumerate}
\end{theorem}
\begin{proof} Observe that (1) implies (2) by Theorem~\ref{t:hlocfac} (2),
(2) implies (3) by Theorem~\ref{t:sfp2} (1), and (3) implies (4)
is trivial.  Assume that $R$ is a Pr\"ufer domain with the
property that each locally divisorial ideal is divisorial. Then it
has the radical trace property by Lemma~\ref{l:locdivrtp}.

Now let $P$ be a nonzero nonmaximal branched prime. Since $R$ has
the radical trace property, $P$ is the radical of a finitely
generated ideal $I$ by \cite[Theorem 23]{lu}. If $M$ is a maximal
ideal that contains $P$, then $J=IR_M\bigcap R$ is locally
divisorial by Lemma~\ref{l:mnlocdiv}. Hence by
Theorem~\ref{t:fgcondiv}, $M$ is the unique maximal ideal that
contains $P$.

Since each unbranched prime must contain a nonzero branched prime,
each nonzero prime is contained in a unique maximal ideal. Thus
$R$ is $h$-local by \cite[Proposition 3.4]{o}. \end{proof}

Our next result adds another equivalence to the $h$-local property
for Pr\"ufer domains.

\begin{theorem} \label{t:converse} Let $R$ be a Pr\"ufer domain with finite
character, and suppose that $R$ has the weak factorization
property. Then $R$ is $h$-local.
\end{theorem}
\begin{proof}  We shall make frequent use of the fact, which
follows easily from \cite[Theorem 1]{gh}, that a Pr\"ufer domain
with finite character satisfies both the \#- and \#\#-properties.
To show that $R$ is $h$-local, it suffices to show that each
nonzero prime ideal is contained in a unique maximal ideal.
Suppose to the contrary that $R$ has a prime ideal $P$ contained
in more than one maximal ideal. Since $R$ has finite character,
$P$ is contained in only finitely many maximal ideals, say $M_1,
\ldots, M_n$, $n>1$. Let $\{P_{\alpha}\}$ denote the set of prime
ideals of $R$ which contain $P$ and are contained in $M_1 \cap
(\bigcup_{j=2}^n M_j)$. This is a chain of prime ideals, and so
$P_1=\bigcup_{\alpha} P_{\alpha}$ is a prime ideal; moreover, $P_1
\sub M_1$, and, by prime avoidance, $P_1 \sub M_i$ for some $i>1$.
One sees easily that $P_1$ is maximal with respect to being
contained in $M_1$ and at least one other maximal ideal. Hence we
may as well assume that $P$ has this property.

Denote by
$\{N_{\alpha}\}$ the set of maximal ideals of $R$ which do not
contain $P$.  Set $T=\bigcap_{j
> 1} R_{M_j} \cap (\bigcap_{\alpha}R_{N_{\alpha}})$.  Since $R$
has finite character, we may find a finitely generated ideal $I$
with the property that $M_1$ is the only maximal ideal containing
$I$. For $x \in I^{-1}$, we have $I \sub (R:_R x)$, so that $(R:_R
x)$ is contained in $M_1$ but no other maximal ideal of $R$.  It
follows that $x \in T$. Hence $I^{-1} \sub T$, and since $I$ is
invertible, $I \supseteq T^{-1}$.  In particular, $M_1 \supseteq
T^{-1}$.

By \cite[Corollary 3.1.8 and Theorem 3.1.2]{fhp}, $P^{-1} = R_P
\cap (\bigcap_{\alpha}R_{N_{\alpha}})$.  In particular $P^{-1}
\supseteq T$.  By Proposition~\ref{p:nonmaxdiv}, $P$ is
divisorial.  Hence $P \sub T^{-1}$.  We claim, in fact, that
$P=T^{-1}$.  Suppose not. Then shrink $M_1$ to a prime ideal $Q$
minimal over $T^{-1}$.  By the maximality property of $P$ and the
fact that $R$ has the \#\#-property, we may choose a finitely
generated ideal $J$ contained in $Q$ such that $M_1$ is the only
maximal ideal of $R$ containing $J$.  As in the preceding
paragraph, we have $T^{-1} \sub J$.  In fact, $T^{-1} \sub J^n$
for each positive integer $n$.  Hence in $R_{M_1}$, we have that
$T^{-1}R_{M_1}$ is contained in the prime ideal $\bigcap_{n \ge 1}
J^nR_{M_1}$ of $R_{M_1}$.  This prime ideal is $Q_0R_{M_1}$ for
some prime ideal $Q_0$ of $R$, and we must have $P \sub T^{-1}
\sub Q_0 \subsetneqq Q$, a contradiction. Thus $P=T^{-1}$, as
claimed.

We next claim that $T$ is a fractional ideal of $R$ which is not
divisorial. Otherwise, the fact that $P=T^{-1}$ implies that
$P^{-1}=T$. However, observe that $T \sub R_{M_2}$, and so it
suffices to show that $P^{-1} \nsubseteq R_{M_2}$.  To see this,
observe by the $\#$-property, $R_{M_1} \cap
(\bigcap_{\alpha}R_{N_{\alpha}}) \nsubseteq R_{M_2}$.  Since
$P^{-1}=R_P \cap (\bigcap_{\alpha}R_{N_{\alpha}}) \supseteq
R_{M_1} \cap (\bigcap_{\alpha}R_{N_{\alpha}})$, we also have
$P^{-1} \nsubseteq R_{M_2}$.  Thus $T$ is not divisorial. Note
that $P^{-1}=T^v \ne T$.

Now consider a possible factorization of $T$: $T=T^v \cdot \Pi$,
where $\Pi$ is a product of maximal ideals.  Then $T=P^{-1} \Pi$.
Since $P^{-1} \sub R_{N_{\alpha}}$, we have $N_{\alpha}P^{-1} \ne
P^{-1}$ (note that $P^{-1}$ is a ring).  If $N_{\alpha}$ appears
as part of $\Pi$, then $1 \in T = P^{-1} \Pi \sub
P^{-1}N_{\alpha}$, a contradiction.  Hence no $N_{\alpha}$ appears
in $\Pi$.  On the other hand, we claim that $M_iP^{-1}=P^{-1}$ for
each $i$. Otherwise, $P^{-1}$ contains a prime ideal $L$
contracting to $M_i$ in $R$, from which it follows that the
valuation domains $(P^{-1})_L$ and $R_{M_i}$ must coincide.
However, the argument in the preceding paragraph showing that
$P^{-1} \nsubseteq R_{M_2}$ can easily be adapted to show that
$P^{-1} \nsubseteq R_{M_i}$.  Hence the claim is true, and we have
$T=P^{-1}\Pi =P^{-1}$, a contradiction. This completes the proof.
\end{proof}

The situation with respect to the weak factorization property is
dramatically different.  Suppose that $R$ is an almost Dedekind
domain with exactly one nondivisorial maximal ideal--see
\cite[Example 42.6]{g1}.  Then $R$ is certainly not $h$-local, but
Theorem~\ref{t:adfact} below implies that $R$ has the weak
factorization property.

\begin{lemma} \label{l:divloc} Let $R$ be an almost Dedekind
domain, let $P$ be an invertible maximal ideal of $R$, and let $I$
be a nonzero ideal of $R$.  Then $I^vR_P=IR_P$.
\end{lemma}
\begin{proof} Since $P$ is invertible, so is $P^i$ for each
$i=1,2, \ldots .$  Hence $I \sub P^i$ if and only if $I^v \sub
P^i$. Since $R_P$ is a rank one discrete valuation domain, we have
$IR_P=P^nR_P$ for some $n \ge 0$.  Since $P^n$ is primary, we then
have $I \sub IR_P \cap R \sub P^nR_P \cap R=P^n$.  Note that $I
\nsubseteq P^{n+1}$.  It follows that $I^v \sub P^n$ and hence
that $I^vR_P=P^nR_P=IR_P$.
\end{proof}

\begin{theorem} \label{t:adfact} Let $R$ be an almost
Dedekind domain, and let $I$ be a nonzero ideal of $R$ which is
contained in only finitely many nondivisorial maximal ideals of
$R$.  Then $I=I^v \cdot \Pi$, where $\Pi$ is a product of maximal
ideals.  Thus, if $R$ is an almost Dedekind domain in which each
nonzero ideal is contained in only finitely many nondivisorial
maximal ideals, then $R$ has the weak factorization property.
\end{theorem}
\begin{proof} Denote by $M_1, \ldots, M_n$ the non-invertible
maximal ideals which contain $I$.  For $M \in \{M_i\}$, we have
$IR_M=M^rR_M$ and $I^vR_M=M^sR_M$ for integers $r,s$ with $0 \le s
\le r$.  Hence $IR_M=I^vM^{r-s}R_M$.  Therefore, for each $i=1,
\ldots, n$, we have a nonnegative integer $t_i$ with
$IR_{M_i}=I^vM^{t_i}R_{M_i}$.  We claim that $I=I^v \cdot
\prod_{i=1}^nM_i^{t_i}$.  We verify this locally.  Let $P$ be a
maximal ideal of $R$.  If $P=M_j$ for some $j$, then
$$IR_P=IR_{M_j}=I^vM_j^{t_j}R_{M_j}=I^v \cdot
(\prod_{i=1}^nM_i^{t_i}R_{M_j}) = (I^v \cdot
\prod_{i=1}^nM_i^{t_i})R_P.$$ If $P \notin \{M_i\}$ and $P$ is
invertible, then, applying Lemma~\ref{l:divloc}, we have
$$IR_P=I^vR_P=(I^v \cdot \prod_{i=1}^n M_i^{t_i})R_P.$$  Finally, if
$P \notin \{M_i\}$ and $P$ is non-invertible, then $I \nsubseteq
P$, so that
$$IR_P=R_P=I^vR_P=(I^v \cdot \prod_{i=1}^n M_i^{t_i})R_P.$$
\end{proof}

Thus any almost Dedekind domain with only finitely many
nondivisorial maximal ideals has the weak factorization property
by Theorem~\ref{t:adfact}.  In fact, it is possible to give
examples of almost Dedekind domains which have infinitely many
nondivisorial maximal ideals but in which each nonzero ideal is
nonetheless contained in only finitely many nondivisorial maximal
ideals--see Example~\ref{e:infnondiv} below.

The next result shows that the integers $t_i$ in the proof of
Theorem~\ref{t:adfact} cannot be ``controlled''.

\begin{proposition} \label{p:adfactarb} Let $R$ be an almost
Dedekind domain, let $M_1 \ldots, M_n$ be distinct non-invertible
maximal ideals of $R$, and let $r_1, \ldots, r_n, s_1, \ldots,
s_n$ be integers with $0 \le s_i \le r_i$.  Then there is a
nonzero ideal $I$ of $R$ such that $I=I^v \cdot \prod_{i=1}^n
M_i^{r_i-s_i}$, and for each $j$, $IR_{M_j}=M_j^{r_j}R_{M_j}$ and
$I^vR_{M_j}=M_j^{s_j}R_{M_j}$.
\end{proposition}
\begin{proof} Note that $M_i \ne M_i^2$ for each $i$ (since this
is true locally).  Hence by ``extended'' prime avoidance
\cite[Theorem 81]{k}, we may pick $a_i \in M_i \sm (\bigcup_{j \ne
i}M_j \cup M_i^2)$. Note that we then have
$M_iR_{M_i}=a_iR_{M_i}$.  Set $I=\prod_{i=1}^n
a_i^{s_i}M_i^{r_i-s_i}$. Since the $M_i$ are non-divisorial, we
have $I^v=\prod_{i=1}^n a_i^{s_i}R$ and hence $I=I^v \cdot
\prod_{i=1}^n M_i^{r_i-s_i}$. Moreover, for each $j$,
$IR_{M_j}=a_j^{s_j}M_j^{r_j-s_j}R_{M_j}=M_j^{r_j}R_{M_j}$, and
$I^vR_{M_j}=a_j^{s_j}R_{M_j}=M_j^{s_j}R_{M_j}$.
\end{proof}

\bigskip


\section{Effects of the strong factorization property}

       Let $D$ be an integral domain with quotient field $K$. Let
$\boldsymbol{\overline{F}}(D)$ denote the set of all nonzero
$D$--submodules of $K$, and let $\boldsymbol{F}(D)$ be the set of
all nonzero fractional ideals of $D$, i.e., $E \in
\boldsymbol{F}(D)$ if $E \in \boldsymbol{\overline{F}}(D)$ and
there exists a nonzero $d \in D$ with $dE \subseteq D$. Let
$\boldsymbol{f}(D)$ be the set of all nonzero finitely generated
$D$--submodules of $K$. Then, obviously $\boldsymbol{f}(D)
\subseteq \boldsymbol{F}(D) \subseteq
\boldsymbol{\overline{F}}(D)$.  A \emph{semistar operation on $D$}
is a map $*: \boldsymbol{F}(D) \to \boldsymbol{F}(D)$, such that,
for each nonzero element $x \in K$ and for each $E,F \in
\boldsymbol{F}(D)$, we have:
\begin{enumerate}
\item $(xE)=(xE)^*$,
\item $E^* \sub F^*$ whenever $E \sub F$, and
\item $E \sub E^*$ and $(E^*)^*=E^*$.
\end{enumerate}

     The semistar operation $*$ on $D$ is called a \it
     (semi)star operation on \rm $D$ if $D^\ast = D$.  (The use of
     the term ``(semi)star'' is due to the fact that, when $D=D^*$, $*$
     is not really a star operation since it remains defined on the
     $D$-submodules of $K$ and not only on the fractional ideals.)

A \emph{localizing system on $D$} is a set $\mc F$ of ideals of
$D$ such that:
\begin{enumerate}
\item if $I \in \mc F$ and $J$ is an ideal of $D$ with $I \sub J$,
then $J \in \mc F$, and
\item if $I \in \mc F$ and $J$ is an ideal of $D$ with $(J:_D a)
\in \mc F$ for each $a \in I$, then $J \in \mc F$.
\end{enumerate}

It is easily seen that a localizing system $\mc F$ is a
multiplicative system of ideals and that $D_{\mc F} := \{x \in K
\mid xI \sub D \text{ for some } I \in \mc F\}$ is an overring of
$D$.  For background on localizing systems, see \cite{fh}, and for
background on semistar operations, see \cite{om} and \cite{fh}.

Now set
       $$ \begin{array}{rl}

       {\mathcal F}^v :=& \hskip -4pt \{ I\mid  I \mbox{ ideal of } D\,,\; I^v =
       D \} \,,\\
       \Pi^v :=& \hskip -4pt   \{ Q \in \Spec(D) \mid Q^v \neq D \mbox{ and } Q
            \neq 0 \}\,,\\

            {\mathcal F}(\Pi^v) :=& \hskip -4pt  \{ I\mid  I \mbox{ ideal of } D\,,\;
            I \not\subseteq Q\,,\; \mbox{ for each } Q \in  \Pi^v\}\,.
            \end{array}
     $$

\begin{lemma} \label{l:1}
           \bara
           \bf \item\it  ${\mathcal F}^v$ is a localizing
           system of $D$ (called \rm the localizing system associated
           to the $v$--operation\it ).

           \bf \item\it The operation $\overline{v} :=
           \ast_{{\mathcal F}^{\!v} }$ defined, for each $E \in
           \overline{\boldsymbol{F}}(D)$, as follows:
           $$
           E^{\overline{v}} := \bigcup \{(E:I) \mid I \in {\mathcal
           F}^v\} \,,$$
           is a (semi)star operation defined on $D$ which is stable
           (i.e. $(E\cap F)^{\overline v} = E^{\overline v}\cap F^{\overline
           v}$, for all $E, F \in  \overline{\boldsymbol{F}}(D)$), and
           it is the largest stable (semi)star operation on $D$.

           \bf \item \it The operation ${v}_{sp} :=
           \ast_{\Pi^v}$ defined, for each $E \in
           \overline{\boldsymbol{F}}(D)$, as follows:
           $$
           E^{{v}_{\!s\!p} } := \bigcap \{ED_{Q}\mid Q \in
           \Pi^v \} \,,$$
           is a semistar operation defined on $D$ (called \rm the
           spectral semistar operation associated to the
           $v$--operation\it ) and  ${\overline v} \leq {v}_{sp} $.

           \bf \item \it
           $$
          {\mathcal F}^{{v}_{\!s\!p} } :=\{ I\mid  I \mbox{ ideal of } D\,,\; I^{{v}_{\!s\!p} } =
       D \} = {\mathcal F}(\Pi^v) \,.
       $$

         \bf \item \it  The following are equivalent:

         \begin{enumerate}
         \bf \item[(i)] \it \ ${v}_{sp}$ is a (semi)star operation on $D$\,;

         \bf \item[(ii)] \it \ ${v}_{sp} \leq v$\,;

         \bf \item[(iii)] \it \  $D = \bigcap \{D_{Q} \mid Q \in \Pi^v \}$\,.
         \end{enumerate}
         \eara

\end{lemma}
\begin{proof} Statements (1), (2), and (3) follow from
  \cite[Proposition~2.8, Theorem~2.10 (B), Proposition~
       3.7 (1), and Proposition~4.11 (2)]{fh}.  Statements (4) and
       (5) are easy consequences of the definitions.
  \end{proof}

\begin{remark} \label{r:rem2}
\rm Note, with respect to Lemma \ref{l:1}
            (2),  that $\overline{v} \leq v$ and so
            $D^{\overline{v}}  = D^v = D$, hence $\overline{v} $ is a
            (semi)star operation on $D$.  As a matter of fact, if $x \in
            E^{\overline{v}} =$ \ $ \bigcup \{(E:I) \mid I \in {\mathcal
           F}^v\}$ then,  for some $I \in {\mathcal
           F}^v$, we have that $I \subseteq (E:_{D}xD)$, thus  $(E:_{D}xD)\in
            {\mathcal
           F}^v$. Therefore, $ D = (E:_{D}xD)^v \subseteq (E^v
           :_{D}xD)$, hence necessarily $1 \in (E^v
           :_{D}xD)$, thus $x \in E^v$.
\end{remark}

\begin{proposition} \label{p:pr2} Assume that $D$ is an $h$-local Pr\"ufer
domain. Then:
           \bara

              \bf \item \it  $\overline{v} = v$\,.

           \bf \item \it
            The following statements are
            equivalent:

            \begin{enumerate}
            \bf \item[(i)] \it The $v$--operation is quasi--spectral (i.e.
            for each nonzero ideal $I$ of $D$, with $I^v \neq D$,
            there exists a prime ideal $Q$ of $D$ such that $I
            \subseteq Q$ and $Q = Q^v$);

            \bf \item[(ii)] \it \ ${v}_{sp} \leq v$\,;

           \bf \item[(iii)] \it  \  $D = \bigcap \{D_{Q} \mid Q \in \Spec(D)\,,
           \; Q = Q^v\}$\,;

            \bf \item[(iv)] \it  \ ${\overline v} ={v}_{sp} =v$\,;

            \bf \item[(v)] \it  \  ${\mathcal F}^v = {\mathcal F}(\Pi^v)$\,.
            \end{enumerate}
           \eara
\end{proposition}
           \begin{proof} (1) It is easy to see that $M^{\ol v}=R$ for
each nondivisorial maximal ideal $M$.  Hence if $I$ is a nonzero
ideal of $R$, the factorization $I=I^vM_1 \cdots M_n$ yields
$I^{\ol v}=(I^v M_1 \cdots M_n)^{\ol v}=(I^v)^{\ol v}=I^v$.

           (2)  These equivalences follow from (1), Theorem \ref{t:sfp2},
           Lemma~\ref{l:1} (5), and
           \cite[Proposition 4.8 and Theorem 4.12 (2)]{fh}.
            \end{proof}

\begin{remark}
           If $V$ is a valuation domain whose maximal ideal $N$
           is idempotent but branched, then $V$ does not satisfy
           any of the (equivalent) conditions in Proposition
           \ref{p:pr2} (2). On the other hand, if $D$ is an $h$-local
           Pr\"ufer domain with non-idempotent maximal ideals, then
           each nonzero ideal of $D$ is divisorial \cite[Theorem 5.1]{h};
           in this case, the (equivalent) conditions in Proposition
           \ref{p:pr2} (2) hold trivially.
\end{remark}


We next study how factorization of an ideal $I$ affects the
factorization of its radical and how factorization of ideals $I$
and $J$ affect the factorization of $IJ$, $I \cap J$, and $I+J$.

\begin{proposition} \label{p:prod}  Let $R$ be an $h$-local
Pr\"ufer domain, and let $I,J$ be nonzero ideals of $R$. Suppose
that $I,J$ have the following factorizations as in
Definition~\ref{d:sfp}:
\begin{align} I &= I^vM_1 \cdots M_k M_{k+1} \cdots M_m H_1
\cdots H_r \text{\quad and} \notag\\
J &= J^vN_1 \cdots N_l N_{l+1} \cdots N_n H_1 \cdots H_r, \notag
\end{align} where the $H_i$ are the nondivisorial maximal ideals
which contain $I+J$ and for which both $IR_{H_i}$ and $JR_{H_i}$
are nondivisorial, $JR_{M_i}$ is principal (including the
possibility that $JR_{M_i}=R_{M_i}$) for $i=1, \ldots, k$, and
divisorial but not principal for $i=k+1, \ldots, m$, and
$IR_{N_i}$ is principal for $i=1, \ldots, l$ and divisorial but
not principal for $i=l+1, \ldots, n$. Further assume that $P_1,
\ldots, P_u$ are the nondivisorial maximal ideals for which
$IR_{P_i}$ and $JR_{P_i}$ are divisorial but $IJR_{P_i}$ is not
divisorial for each $i$. Then the canonical factorizations of $IJ$
and $I^vJ^v$ are as follows:

\begin{align} IJ &= (IJ)^vM_1 \cdots M_kN_1 \cdots N_l H_1
\cdots H_rP_1 \cdots P_u\\
I^vJ^v &= (IJ)^vP_1 \cdots P_u.
\end{align}
\end{proposition}
\begin{proof} (1)  For each $i=1, \ldots, k$, we have elements
$x_i,y_i \in R$ with $IR_{M_i} = x_iM_iR_{M_i}$ and
$JR_{M_i}=y_iR_{M_i}$, so that $IJR_{M_i}=x_iy_iM_iR_{M_i}$.  Thus
$IJR_{M_i}$ is not divisorial, and each of these $M_i$ must appear
in the factorization of $IJ$.  Similarly, $N_1, \ldots, N_l$ must
appear.  For $i=k+1, \ldots, m$, there is an element $z_i \in R$
with $IJR_{M_i}=z_iM_iJR_{M_i}=z_iJR_{M_i}$; the second equality
follows from the fact that in a valuation domain with maximal
ideal $Q$ a nonprincipal ideal $K$ satisfies $K=KQ$ (see
Remark~\ref{r:sfp}). In this case, $IJR_{M_i}$ is divisorial, and
so $M_i$ does not appear in the factorization of $IJ$. Similarly,
$N_{l+1}, \ldots, N_n$ do not appear.  For $H \in
\{H_i\}_{i=1}^r$, since both $IR_H$ and $JR_H$ are nondivisorial,
there are elements $x,y$ with $IJR_H=xHyHR_H=xyHR_H$ (note that
$H$ is idempotent by Theorem~\ref{t:sfp2} (2)); this is not
divisorial, so each $H_i$ must appear. Finally, it is clear that
the $P_i$ must appear and that no other maximal ideals can appear.

(2)  First, observe that if $Q$ is a nondivisorial maximal ideal
for which $IR_Q$, $JR_Q$, and $IJR_Q$ are all divisorial, then by
\cite[Lemma 2.3]{bs}, $I^vJ^vR_Q=(IR_Q)^v(JR_Q)^v=IJR_Q$, which is
divisorial. Hence no such $Q$ appears in the factorization of
$I^vJ^v$.  Let $M \in \{M_i\}_{i=1}^m$.  Then there is an element
$x \in R$ with $I^vJ^vR_M=I^vJR_M=J(IR_M)^v=J(xMR_M)^v=JxR_M$,
which is divisorial.  Thus no $M_i$ appears; similarly, no $N_i$
appears. For $H \in \{H_i\}_{i=1}^r$, we have an element $y \in R$
with $I^vJ^vR_H=I^v(JR_H)^v=I^v(yMR_H)^v=I^vyR_H$, which is
divisorial. Thus no $H_i$ appears.
\end{proof}

\begin{proposition} \label{p:intersection} Let $R$ be an $h$-local Pr\"ufer
domain, and let $I,J$ be nonzero ideals of $R$. Suppose that $I,J$
have the following factorizations as in Definition~\ref{d:sfp}:
\begin{align} I &= I^vM_1 \cdots M_k M_{k+1} \cdots M_m H_1
\cdots H_r \text{\quad and} \notag\\
J &= J^vN_1 \cdots N_l N_{l+1} \cdots N_n H_1 \cdots H_r, \notag
\end{align} where the $H_i$ are the nondivisorial maximal ideals which contain
$I+J$ and  for which both $IR_{H_i}$ and $JR_{H_i}$ are
nondivisorial, $IR_{M_i} \sub JR_{M_i}$ for $i=1, \ldots, k$,
$IR_{M_i} \nsubseteq JR_{M_i}$ for $i=k+1, \ldots, m$, $JR_{N_i}
\sub IR_{N_i}$ for $i=1, \ldots, l$, and $JR_{N_i} \nsubseteq
IR_{N_i}$ for $i=l+1, \ldots, n$.  Then $I \cap J$ has the
following factorization $$I \cap J=(I \cap J)^vM_1 \cdots M_kN_1
\cdots N_lH_1 \cdots H_r.$$
\end{proposition}
\begin{proof} For $i=1, \ldots, k$, $(I \cap J)R_{M_i}=IR_{M_i}$,
and so $M_i$ appears in the factorization of $I \cap J$. Moreover,
for $j > k$, $(I \cap J)R_{M_j}=JR_{M_j}$; since $JR_{M_j}$ is
divisorial, $M_j$ does not appear.  The $N_i$ are handled
similarly.  Finally, it is straightforward to show that the $H_i$
appear and that no other maximal ideals can appear.
\end{proof}

\begin{proposition} \label{p:sum} Let $R$ be an $h$-local Pr\"ufer
domain, and let $I,J$ be nonzero ideals of $R$. Then:
\begin{enumerate}
\item If $I$ and $J$ are divisorial, then $I+J$ is divisorial.
\item In general, $(I+J)^v=I^v+J^v$.
\item Let $I$ and $J$ have the following factorizations as in Definition~\ref{d:sfp}:
\begin{align} I &= I^vM_1 \cdots M_k M_{k+1} \cdots M_m H_1 \cdots H_r \notag\\
    J &= J^vN_1 \cdots N_l N_{l+1} \cdots N_n H_1 \cdots H_r, \notag
\end{align}
where the $H_i$ are the nondivisorial maximal ideals which contain
$I+J$ and  for which both $IR_{H_i}$ and $JR_{H_i}$ are
nondivisorial, $IR_{M_i} \sub JR_{M_i}$ for $i=1, \ldots, k$,
$IR_{M_i} \nsubseteq JR_{M_i}$ for $i=k+1, \ldots, m$, $JR_{N_i}
\sub IR_{N_i}$ for $i=1, \ldots, l$, and $JR_{N_i} \nsubseteq
IR_{N_i}$ for $i=l+1, \ldots, n$.  Then the factorization of $I+J$
is $$I+J = (I+J)^vM_{k+1} \cdots M_mN_{l+1} \cdots N_nH_1 \cdots
H_r.$$
\end{enumerate}
\end{proposition}

\begin{proof} (1) Let $M$ be a maximal ideal of $R$.  By Theorem \ref{t:sfp2},
both $IR_M$ and $JR_M$ are divisorial.  Since $(I+J)R_M$ is equal
to one of these, it is divisorial.  Hence $I+J$ is divisorial,
again by Theorem \ref{t:sfp2}.

(2) Using (1), we have $(I+J)^v = (I^v + J^v)^v = I^v + J^v$.

(3)  Using the same reasoning as in the proof of
Proposition~\ref{p:prod} (1), we see easily that each $H_i$ must
appear in the factorization of $I+J$. Similarly, for $M \in
\{M_i\}_{i=k+1}^m$, we have $(I+J)R_M = IR_M$, so these $M_i$ must
appear.  Each $N_i$, $i=l+1, \ldots, n$, must also appear.  The
same reasoning shows that none of the other $M_i$ or $N_j$ can
appear, and it is clear that no other maximal ideals can appear.
\end{proof}

\begin{proposition} \label{p:rad} Let $R$ be an $h$-local
Pr\"ufer domain, and let $I$ be an ideal of $R$ with factorization
(as in Definition~\ref{d:sfp}) $$I=I^vM_1 \cdots M_l M_{l+1}
\cdots M_k M_{k+1} \cdots M_n,$$ where $M_1, \ldots, M_k$ are
minimal over $I$, $I^v \sub M_i$ for $i=1, \ldots, l$, $I^v
\nsubseteq M_i$ for $i=l+1, \ldots, k$, and $M_i$ is not minimal
over $I$ for $i=k+1, \ldots, n$.  Let $\{N_1, \ldots, N_r\}$
denote the (possibly empty) set of nondivisorial maximal ideals
that are minimal over $I$ and are such that $IR_{N_i}$ is
divisorial. Then:
\begin{enumerate}
\item The factorization of $\rad I$ is $\rad I = (\rad I)^v M_1 \cdots M_k N_1 \cdots N_r$.
\item $(\rad I)^v = (\rad I^v)^v$.
\item The factorization of $\rad I^v$ is $\rad I^v = (\rad I)^v M_1
\cdots M_lN_1 \cdots N_r$. \end{enumerate}
\end{proposition}
\begin{proof} (1) For $i=1, \ldots, k$, $(\rad I)R_{M_i}=
M_iR_{M_i}$, so $M_i$ must appear in the factorization of $\rad
I$.  Also, since $(\rad I)R_{N_i}=N_iR_{N_i}$, each $N_i$ must
appear. For any other nondivisorial maximal ideal $P$ containing
$I$, $P$ is not minimal over $I$, whence $(\rad I) R_P$ is a
nonmaximal, and hence divisorial, prime ideal of $R_P$.

(3)  We have $(\rad I)^v = (\rad (I^v \prod_{i=1}^n M_i))^v= (\rad
I^v \cap \prod_{i=1}^n M_i)^v = (\rad I^v)^v$, with the last
equality following from the fact that the $v$-operation is stable
in the presence of strong factorization (Proposition~\ref{p:pr2}).

(4)  For $Q \in \{M_i\}_{i=1}^l$, it is clear that $Q$ is minimal
over $I^v$.  For $Q \in \{N_i\}_{i=1}^r$, use the fact that $R$ is
$h$-local to obtain $IR_Q=I^vR_Q$. Since $I \sub Q$, we must have
$I^v \sub Q$, and, again, $Q$ is minimal over $I^v$. In either
case, we therefore have $(\rad I^v)R_Q=QR_Q$, which is
nondivisorial, whence $Q$ must appear in the factorization of
$\rad I^v$.  It is clear that no other maximal ideals can appear.
\end{proof}

\begin{proposition} \label{p:trace} Let $R$ be an $h$-local Pr\"ufer
domain.  Let $I$ be a nonzero ideal of $R$, and suppose that the
factorization of $I$ (as in Definition~\ref{d:sfp}) is $I=I^vM_1
\cdots M_n$. Let $P_1, \ldots, P_u$ be the nondivisorial maximal
ideals containing $II^{-1}$ for which $IR_{P_i}$ is divisorial but
$II^{-1}R_{P_i}$ is nondivisorial.  Then the factorization of
$II^{-1}$ is $II^{-1}=(II^{-1})^vM_1 \cdots M_n P_1 \cdots P_u$.
\end{proposition}

\begin{proof} For $M \in \{M_i\}_{i=1}^n$, there is an element $x \in R$
with $$II^{-1}R_M = xMI^{-1}R_M = xM(IR_M)^{-1} = xM(xMR_M)^{-1} =
MR_M,$$ where the second equality follows from the fact that $R$
is $h$-local \cite[Lemma 2.3]{bs}.  Hence each $M_i$ must appear.
It is clear that each $P_i$ must appear and that no other maximal
ideals can appear.
\end{proof}

\medskip

We observe that the $P_i$ in Propositions \ref{p:prod} and
\ref{p:trace} can actually occur--see Example~\ref{e:prodnotdiv}
below.

             \medskip

We end this section with a result which contains more information
related to Propositions~\ref{p:intersection} and \ref{p:trace}

\begin{proposition} \label{p:maybenew} Let $R$ be an $h$-local
Pr\"ufer domain.  If $I$ is a nondivisorial ideal of $R$ with
factorization $I=I^vM_1 \cdots M_n$ (as in
Definition~\ref{d:sfp}), then
\begin{enumerate}
\item for each $i=1, \ldots, n$, $I^vI^{-1} \nsubseteq M_i$, and
$I^vR_{M_i}$ is principal;
\item $II^{-1}=I^vI^{-1} M_1 \cdots M_n$, and for each $i=1,
\ldots, n$ $M_i$ is minimal over $II^{-1}$ and
$II^{-1}R_{M_i}=M_iR_{M_i}$;
\item there is a finitely generated ideal $J \sub I^v$ with
$I+J=I^v$, and, for any such $J$, $(I \cap J)^v=J$; and
\item for each nonzero ideal $B \sub I^v$, $(I \cap B)^v=B^v$.
\end{enumerate}
\end{proposition}
\begin{proof} Let $M \in \{M_i\}$.  Write $I^vMR_M=IR_M=xMR_M$,
where (we may assume) $x \in I^v$.  Then $I^vR_M=(IR_M)^v=xR_M$,
and by \cite[Theorem 3.10]{o} $I^{-1}R_M=(IR_M)^{-1}=x^{-1}R_M$.
It follows that $I^vI^{-1} \nsubseteq M$.  In particular, $I^vR_M$
is invertible, hence principal, in $R_M$, proving (1).

For (2), from what was just proved, we have $II^{-1}=I^vI^{-1}M_1
\cdots M_n$ with $I^vI^{-1}$ and $M_1 \cdots M_n$ comaximal.  Thus
$II^{-1}R_{M_i}=M_iR_{M_i}$, as desired.

Now let $J=(x_1, \ldots, x_n) \sub I^v$ be such that
$x_iR_{M_i}=I^vR_{M_i}$ for each $i$.  Then for $M \in \{M_i\}$,
we have $I^vR_M=JR_M=(I+J)R_M$.  On the other hand, if $N$ is a
maximal ideal with $N \notin \{M_i\}$, then
$IR_N=(IR_N)^v=I^vR_N$, from which it follows easily that
$(I+J)R_N=I^vR_N$.  Therefore, $I+J=I^v$.  Using
Proposition~\ref{p:pr2} (1), we also obtain $(I \cap J)^v=I^v \cap
J^v=J^v=J$, proving (3).  Statement (4) also follows from
Proposition~\ref{p:pr2} (1).
\end{proof}

\bigskip


\section{Examples}

We begin with a lemma which is probably known but for which we
have no convenient reference.

\begin{lemma} \label{l:pid} For any nonempty set of indeterminates
$\mathcal Z=\{Z_\alpha\}$ and any field $F$, the ring $D=\bigcap
F[\mathcal Z]_{(Z_\alpha)}$ is a PID with $\Max(D)=\{Z_\alpha
D\,|\, Z_\alpha\in \mathcal Z\}$.
\end{lemma}

\begin{proof} Let
$f\in F[\mathcal Z]$. If no $Z_\alpha$ divides $f$ in $F[\mathcal  Z]$, then
$f^{-1}$ is in each localization $F[\mathcal Z]_{(Z_\alpha)}$. Thus a reduced
rational expression $g/f$ from the quotient field of $F[\mathcal Z]$ is in $D$
if and only if no $Z_\alpha$ divides $f$. Thus each element of $D$ has the
reduced form $g/f$ where no $Z_\alpha$ divides $f$. Clearly $g/f$ is a unit of
$D$ if and only if no $Z_\alpha$ divides $g$. It follows that each nonzero
prime ideal  of $D$ is principal of the form $Z_\alpha D$ for some (unique)
$Z_\alpha$. Hence $D$ is a PID.
\end{proof}

\begin{example} \label{e:infnondiv}An example of an almost Dedekind domain
$D$ with infinitely many nondivisorial maximal ideals such that
$D$ has the weak factorization property.

Notation:
\begin{enumerate}
\item For each $n \ge 1$, let $\mathcal X_n=\prod_{i>0} X_{n,i}$  where   $\{X_{n,i}\,|\, 1\le i, \, 1\le n\}$  is a countably infinite set of  algebraically independent indeterminates.  \tpp
\item For each $n$ and each $k\ge 0$, let $\mathcal
X_{n,k}=\prod_{i>k}X_{n,i}$ (so $\mathcal X_{n,0}=\mathcal
X_n$).\tpp
\item Let $E_0=K[\{\mathcal X_n\,|\, 1\le n\}]$   and for each $n$, let
$Q_{n,0}=(\mathcal X_n)E_0$. \tpp
\item For each $k\ge 1$, let $E_k=K[\{X_{n,j} \,|\, 1\le j\le k, 1\le n\}, \{\mathcal X_{n,k} \,|\, 1\le n\}]$, $P_{n,j}=(X_{n,j})E_k$ for $j\le k$ and $Q_{n,k}=(\mathcal X_{n,k})E_k$. \tpp
\item Let $D_0=\bigcap (E_0)_{Q_{n,0}}$ and for $k\ge 1$, let $D_k=\left(\bigcap
(E_k)_{Q_{n,k}}\right) \cap \left(\bigcap (E_k)_{P_{n,j}}\right)$.\tpp
\item Finally let $D=\bigcup D_k$.\end{enumerate}

Then \begin{enumerate}
\item $D$ is an almost Dedekind domain which is also a Bezout domain.
\item Each nonzero ideal is contained in at most finitely many
nondivisorial maximal ideals.
\item $D$ has the weak factorization property.
\end{enumerate}
\vskip8pt

\begin{proof}  Each $D_k$ is a PID. Also it is clear that each maximal ideal
$M$ of $D_k$ contracts to a maximal ideal $N_j$ of $D_j$ for each
$j<k$ and $N_j(D_k)_M=M(D_k)_M$. Moreover, each maximal ideal of
$D_k$ survives in $D_m$ for each $m>k$. Thus by \cite[Theorem
2.10]{ll}, $D$ is an almost Dedekind domain that is also a Bezout
domain -- given a finitely generated ideal $I$ of $D$, $I=I_kD$
where $I_k=I\cap D_k$ for some $k$.

By the proof of \cite[Theorem 2.10]{ll}, each maximal ideal $M$  of $D$ is the union of its contractions to the $D_k$'s. As in  the proof of \cite[Example 3.2]{ll}, $D$ has two distinct types of maximal ideals. For each $X_{n,k}$, the ideal  $M_{n,k}=X_{n,k}D$ is a principal maximal ideal of $D$. The other maximal ideals are those of the form $M_n=\bigcup_{j\ge 0} Q_{n,j}$.   For each $n$, we let $\mathcal F_n=\{M_n,M_{n,1}, M_{n,2},\dots\}$ and call this the family of maximal ideals centered on $\mathcal X_n$.  These are the only maximal ideals of $D$ that contain $\mathcal X_n$ (and each does). Since $D$ is an almost Dedekind domain,  some member of $\mathcal F_n$ is not finitely generated. The only one that is not principal is $M_n$. Thus $M_n$ is not divisorial.

For a nonzero proper ideal $I$, recall that $\Max(R,I)$ is the set
of maximal ideals of $D$ that contain $I$; let us refer to this as
the \emph{support} of $I$. We will show that $\Max(R,I)$ is
contained in a finite union of families $\mathcal F_n$. To this
end, let $f$ be a nonzero nonunit of $D$ and let $D_k$ be the
smallest member of the chain that contains $f$. By the argument
above, $f=ug/v$ with $u$ and $v$ units of $D_k$ and $g$  a finite
product of monomials of the form $\mathcal X_{n,k}$ and $X_{m,i}$
with $i\le k$. Since $u$ and $v$ are units of $D$, the monomials
in $g$ completely determine the families that contain the support
of $f$. Thus $\Max(R,(f))$ is contained in the union of finitely
many families $\mathcal F_n$. Hence the same is true for the
support of each nonzero proper ideal. Moreover, since each family
contains exactly one nondivisorial ideal, each nonzero proper
ideal is contained in at most finitely many nondivisorial maximal
ideals. Therefore, $D$ has the weak factorization property by
Therorem \ref{t:adfact}.
\end{proof}
\end{example}

\bpb

\begin{example} \label{e:sumnotdiv} An example of a Pr\"ufer domain $R$ with the weak factorization property such that $R$ contains ideals
$I,J$ with $I$ and $J$ divisorial but $I+J$ not divisorial.

  We recall the construction of the domain in \cite[Example 3.2]{ll}.

Let $\mc X=\prod_{i > 0} X_i$, where the $X_i$ are indeterminates.
Let $K$ be a field, and for each $n$, let $\mc X_n = \prod_{k \ge
n} X_k$ and $E_n = K[X_1, \ldots, X_{n-1}, \mc X_n]$ ($E_0 = K[\mc
X]$).  Set $P_{n,k}=X_kE_n$, $P_n=\mc X_nE_n$, and
$D_n=\left(\bigcap_{k < n}(E_n)_{P{n,k}}\right) \cap (E_n)_{P_n}$.
Let $Q_{n,k} = P_{n,k}D_n$ and $Q_n = P_nD_n$.  Then each $D_n$ is
a semilocal PID, and $D=\bigcup D_n$ is an almost Dedekind domain with a unique noninvertible maximal ideal.  We also have the following.

\begin{enumerate}
\item $D$ has only  countably many maximal ideals
$M, M_1, M_2, \ldots$, where $M=\bigcup {Q_n}$, and $M_n = X_nD$.
Also $D$ has nonzero Jacobson radical, since $\mc X$ is in each maximal
ideal.  The maximal ideal $M$ is nondivisorial, while the $M_n$'s
are all principal.
\item The ideals $I = \bigcap_{k \ge 1}M_{2k}$ and $J = \bigcap_{k \ge 1}M_{2k-1}$ are (nonzero) divisorial ideals, but $I+J$ is nondivisorial.
\item $D$ has the weak factorization property.
\end{enumerate}

\begin{proof}  Statement (1) is from  \cite[Example 3.2]{ll}.

Since $D$ has nonzero Jacobson radical, $I$ and $J$ are nonzero;
they are divisorial since each $M_n$ is divisorial.  We have $I+J
\sub M$ since each element of $D$ which is contained in infinitely
many $M_n$ is also in $M$ (see either Lemma 2.2 or Theorem 2.5 of \cite{ll}).  In fact, we claim that
$I+J = M$.  Observe that $\mc XR_M = MR_M$ so that $(I+J)R_M =
MR_M$.  Moreover, for each positive even integer $k$ the element
$X_2X_4 \cdots X_k \mc X_{k+1}$ is in $I$ but is a unit in
$D_{M_{k-1}}$; hence $(I+J)D_{M_{k-1}} = ID_{M_{k-1}} =
D_{M_{k-1}} = MD_{M_{k-1}}$.  Applying the same argument to $J$,
we obtain $(I+J)D_{M_k} = MD_{M_{k}}$.  It follows that $I+J = M$,
so $I+J$ is not divisorial.
\end{proof}
 \end{example}

 \begin{example} \label{e:prodnotdiv} An example of a valuation containing
$V$ containing a divisorial $I$ for which $II^{-1}$ is not
divisorial (thus the product of divisorial ideals need not be
divisorial).

Let $(V,M)$ be an valuation domain with value group the additive
rational numbers. Note that $M$ is not principal and therefore not
divisorial.  Let $I$ denote the ideal consisting of those elements
of $V$ having value greater than $\sqrt 2$.  For each positive
rational number $\alpha$, let $x_{\alpha}$ denote an element of
$V$ with value $\alpha$.  Then $I=\bigcap_{\alpha < \sqrt 2}
(x_{\alpha})$. Hence $I$ is divisorial.   However, $I$ is not
(principal hence not) invertible, whence by \cite[Proposition
4.2.1]{fhp} $II^{-1}$ must be a prime ideal of $V$.  Since $V$ is
one-dimensional, we must therefore have $II^{-1}=M$, which is not
divisorial.
 \end{example}

\begin{example} \label{e:nofac} An example of a one-dimensional Bezout
domain which does not have the weak factorization property.

Let $\mathcal X=\prod_{k\ge 0} X_k^{2^k}$ where $\{X_k\}$ is a
countably infinite set of indeterminates. Let $K$ be a field, and
for each integer $n$, let $\mathcal X_n=\prod_{k\ge n}
X_k^{2^{k-n}}$ and $E_n=K[X_0,X_1,\dots, X_{n-1}, \mathcal X_n]$
(with $E_0=K[\mathcal X]$). Let $P_{n,k}=X_kE_n$ for $k<n$,
$P_n=\mathcal X_nE_n$, and $D_n=(\bigcap (E_n)_{P_{n,k}})\cap
(E_n)_{P_n}$. Use $Q_{n,k}$ to denote the extension of $P_{n,k}$
to $D_n$ and $Q_n$ to denote the extension of $P_n$ to $D_n$. Each
$Q_{n,k}$ is principal as is each $Q_n$. Also each $D_n$ is a
semilocal  PID.

Let $D=\bigcup D_n$. Then $D$ is a one-dimensional Bezout domain
with nonzero Jacobson radical. Also, $D$  has countably many
maximal ideals. Of these, all but one is principal. The one that
is not principal is idempotent. This maximal ideal is not the
radical of a finitely generated ideal, so it is non-sharp.  It
follows that from Proposition~\ref{p:nonmaxdiv} (3) that $D$ does
not have the weak factorization property.

\begin{proof}   Let $I=(a_1,a_2,\dots,a_m)$ be a finitely generated proper
ideal of $D$. Let $D_n$ be the smallest ring in $\{D_i\}$ that
contains the set $\{a_1,a_2,\dots,a_m\}$. Since $D_n$ is a PID,
there is an element $a\in I\cap D_n$ such that $I\cap D_n=aD_n$.
In particular, each $a_i$ is in $aD_n$ and it follows that $I=aD$.
Thus $D$ is a Bezout domain.

For  integers $0\le m<n$ and $0\le k<n$, $Q_{n,k}\bigcap D_m=Q_{m,k}$ when
$k<m$ and $Q_{n,k}\bigcap D_m=Q_m$ when $m\le k$. In the first case,
$Q_{m,k}(D_m)_{Q_{n,k}}=Q_{n,k}(D_n)_{Q_{n,k}}$, and in the
second, $Q_m(D_n)_{Q_{n,k}}=Q_{n,k}^j(D_n)_{Q_{n,k}}$ where $j=2^{k-m}$.

Let $f$ be a nonzero member of $D$. Since $D$ is the union of the chain $D_n$
and no nonunit of $D_n$ becomes a unit in a larger $D_m$, $f$ is a nonunit of
$D$ if and only if it is a nonunit in the smallest  $D_n$ that contains it. In
$D_n$, $f$ is a nonunit if and only if has the form $ug/v$ where $u$ and $v$
are polynomials of $E_n$ that are units of $D_n$ and $g$ is a finite
(nonempty) product of the  monomials $X_0$, $X_1$, \dots, $X_{n-1}$ and
$\mathcal X_n$. If the factorization of $g$ does not include a positive power
of $\mathcal X_n$, then for all $m>n$, $f\notin Q_m$. On the other hand, if
the factorization of $g$ does include a positive power of $\mathcal X_n$, then
$f\in Q_m$ for all $m\ge n$. In the latter case, we also have that $f\in
Q_{m,k}$ for all $m>k\ge n$ since $\mathcal X_n=\mathcal
X_m^{2^{m-n}}\,\prod_{k=n}^{m-1} X_k^{2^{k-n}}$

For each $n$ the ideal $M_n=X_nD$ is a height one maximal ideal of
$D$, being the union of the chain of primes $Q_0\subset \cdots
\subset Q_{n-1}\subset Q_{n,n}\subset Q_{n+1,n}\subset \cdots$.
The only other maximal ideal of $D$ is the ideal $M=\bigcup Q_n$,
the union of the chain $\{Q_n\}_{0\le n}$. The height of  $M$ is
also one, so $D$ is one-dimensional. Let $f$ be a nonzero member
of $M$. Then there is an integer $n$ such that $f$ is in  $Q_m$
for each $m\ge n$. But this implies that $f\in Q_{m,k}$ for each
pair $m>k\ge n$.  Since $D$ is a Bezout domain, $M$ cannot be the
radical of a finitely generated ideal.
\end{proof}
\end{example}

\begin{remark} \label{r:divnotloc} It is perhaps worth noting that
the preceding provides an example of a divisorial ideal $J$ in a
Pr\"ufer domain such that $JR_M$ is not divisorial for some
maximal ideal $M$. With the notation above, let $J$ be the
intersection of the principal maximal ideals.  Then $J$ is nonzero
and divisorial. We must have $J \sub M$. Otherwise, for $x \in J
\sm M$ we would have $(M,x)=R$.  However, writing $1=m+rx$, $m \in
M$, $r \in R$ then yields that $M$ is the only maximal ideal
containing $m$, a contradiction.  Since $J$ is a radical ideal, we
must then have $JR_M=MR_M$, which is nondivisorial.
\end{remark}

           \end{document}